\documentclass[12pt,a4paper,reqno]{amsart}

\usepackage{amsmath,amsthm}
\usepackage{amssymb}

\numberwithin{equation}{section}

\theoremstyle{theorem}

\theoremstyle{definition}

\newtheorem{remark}{Remark}

\theoremstyle{remark}

\usepackage{url}
\usepackage[colorlinks,linkcolor=blue,citecolor=blue]{hyperref}

\title{A note on the coefficients of power sums of arithmetic progressions}

\author[J.L. Cereceda]{Jos\'e Luis Cereceda}
\address{%
        Collado Villalba, 28400 -- Madrid, Spain}
\email{jl.cereceda@movistar.es}

\begin{document}

\begin{abstract}
In this note we show a simple formula for the coefficients of the polynomial associated with the sums of powers of the terms of an arbitrary arithmetic progression. This formula consists of a double sum involving only ordinary binomial coefficients and binomial powers. Arguably, this is the simplest formula that can probably be found for the said coefficients. Furthermore, we give an explicit formula for the Bernoulli polynomials involving the Stirling numbers of the first and second kind.
\end{abstract}

\maketitle

\vspace{-.5cm}

\section{Introduction}

Consider the sum of $k$th powers of the terms of an arithmetic progression with first term $r$ and common difference $m$
\begin{equation*}
S_k^{m,r}(n) = r^k + (m+r)^k + (2m+r)^k + \cdots + ( (n-1)m +r)^k,
\end{equation*}
where $k,m,r,$ and $n$ are assumed to be integer variables with $m,n \geq 1$ and $k,r \geq 0$. $S_k^{m,r}(n)$ turns out to be a polynomial in $n$ of degree $k+1$ without constant term, that is, it can be expressed in the form $S_k^{m,r}(n) = \sum_{t=1}^{k+1} c_{k,t}^{m,r} n^t$ for certain rational coefficients $c_{k,t}^{m,r}$.

Griffiths \cite{grif} obtained the following formula for $c_{k,t}^{m,r}$ (in our notation)
\begin{equation}\label{c1}
c_{k,t}^{m,r} = m^k \sum_{j=t}^{k+1} \sum_{i=j}^{k+1} \frac{1}{i} \left( \frac{r}{m} \right)^{j-t} \binom{j}{t} S_1(i,j) S_2(k,i-1),
\end{equation}
for $t =1,\ldots, k+1$, where $S_1(i,j)$ are the (signed) Stirling numbers of the first kind and $S_2(k,i-1)$ are the (unsigned) Stirling numbers of the second kind. On the other hand, the authors of \cite{bazso} gave the following formula for $c_{k,t}^{m,r}$
\begin{equation}\label{c2}
c_{k,t}^{m,r} = \sum_{j=0}^{k} \frac{m^j W_{m,r}(k,j)}{j+1} S_1(j+1,t),
\end{equation}
where $S_1(j+1,t)$ are the (signed) Stirling numbers of the first kind and $W_{m,r}(k,j)$ are the so-called $r$-Whitney numbers of the second kind, which are defined in \cite{mezo} by
\begin{equation*}
W_{m,r}(k,j) = \frac{1}{m^j j!} \sum_{i=0}^{j} (-1)^{j-i} \binom{j}{i} \big( m i +r \big)^k .
\end{equation*}
Moreover, from \cite[Theorems 5 and 6]{ramirez}, it can be shown that
\begin{equation}\label{c3}
c_{k,t}^{m,r} = \frac{(-1)^{k+1-t}}{(k+1)!} \sum_{j=0}^{k} \sum_{s =t}^{k+1} A_{m,r}(k,k-j) \binom{s}{t}
\big( k-j \big)^{s-t} S_1(k+1,s),
\end{equation}
where $S_1(k+1,s)$ are the (signed) Stirling numbers of the first kind, and
\begin{equation*}
A_{m,r}(k,j) = \sum_{i=0}^{j} (-1)^i \binom{k+1}{i} \big[ m(j-i) + r \big]^{k}.
\end{equation*}
It is to be noted that the polynomial $S_k^{m,r}(n) = \sum_{t=1}^{k+1} c_{k,t}^{m,r} n^t$ with coefficients given by \eqref{c1}, \eqref{c2}, and \eqref{c3} can equivalently be expressed, for $k \geq 0$, as
\begin{align*}
S_k^{m,r}(n) & = m^k \sum_{j=0}^{k} j! S_2(k,j) \left[ \binom{n+\frac{r}{m}}{j+1} - \binom{\frac{r}{m}}{j+1}
\right],  \\
S_k^{m,r}(n) & = \sum_{j=0}^{k} m^j j! W_{m,r}(k,j) \binom{n}{j+1}, \\[-2mm]
\intertext{and}
S_k^{m,r}(n) & = \sum_{j=0}^{k} A_{m,r}(k,k-j) \binom{n+j}{k+1},
\end{align*}
respectively. For other explicit formulas concerning the polynomial $S_k^{m,r}(n)$ as a whole see, for example, \cite{gaut,hirs,tsao}.

\section{A simpler formula for $c_{k,t}^{m,r}$}

The main purpose of this note is to show that a simpler formula for $c_{k,t}^{m,r}$, as compared to the formulas in \eqref{c1}, \eqref{c2}, and \eqref{c3}, can be derived from the following alternative expression for $S_k^{m,r}(n)$ (see, e.g., \cite[Equation 11]{bazso} and \cite[Note 3]{cere})
\begin{equation}\label{ber}
S_k^{m,r}(n) = \frac{m^k}{k+1} \sum_{t=1}^{k+1} \binom{k+1}{t} B_{k+1-t}\left( \frac{r}{m}\right) n^t, \quad k \geq 0,
\end{equation}
where $B_k(x)$ are the Bernoulli polynomials with generating function
\begin{equation*}
\frac{t e^{xt}}{e^t -1} = \sum_{k=0}^{\infty} B_k(x) \frac{t^k}{k!}, \quad ( |t| < 2\pi ).
\end{equation*}
Equation \eqref{ber} generalizes the well-known formula for the ordinary power sum $S_k^{1,0}(n) = 0^k + 1^k + \cdots + (n-1)^k$, namely
\begin{equation*}
S_k^{1,0}(n) = \frac{1}{k+1} \sum_{t=1}^{k+1} \binom{k+1}{t} B_{k+1-t} n^t, \quad k \geq 0,
\end{equation*}
where $B_k = B_k(0)$ are the Bernoulli numbers.

Indeed, by using the following specific representation for $B_k(x)$ (see, e.g, \cite[Equation (38)]{todorov} and \cite[Equation (18)]{sri})
\begin{equation}\label{ber2}
B_{k}(x) = \sum_{j=0}^{k} \sum_{i=0}^{j} \frac{(-1)^{i}}{j+1} \binom{j}{i} \big( i + x \big)^{k}, \quad k \geq 0,
\end{equation}
and setting $x \to \frac{r}{m}$, we quickly obtain from \eqref{ber}
\begin{equation}\label{myc}
c_{k,t}^{m,r} = \frac{m^k}{k+1} \binom{k+1}{t} \sum_{j=0}^{k+1-t} \sum_{i=0}^{j} \frac{(-1)^{i}}{j+1} \binom{j}{i}
\left( i + \frac{r}{m} \right)^{k+1-t},
\end{equation}
which gives us $c_{k,t}^{m,r}$ as a double sum involving only ordinary binomial coefficients and binomial powers. Alternatively, we can write \eqref{myc} as
\begin{equation}\label{myc2}
c_{k,t}^{m,r} = \frac{m^{t-1}}{k+1} \binom{k+1}{t} \sum_{j=0}^{k+1-t} (-1)^j \frac{m^j j!}{j+1} W_{m,r}(k+1-t,j).
\end{equation}
Formally, this last formula is simpler than those in \eqref{c1}, \eqref{c2}, and \eqref{c3} because each of the numbers $S_1(k,j)$, $S_2(k,j)$, and $A_{m,r}(k,j)$, carries at least one extra summation of its own. In particular, the (signed) Stirling numbers of the first kind are given explicitly by
\begin{equation*}
S_1(k,j) = \sum_{t=0}^{k-j} \sum_{r=0}^{t} \frac{(-1)^r}{t!} \binom{k+t-1}{k+t-j} \binom{2k-j}{k-t-j} \binom{t}{r}
r^{k-j+t},
\end{equation*}
so that each of the formulas \eqref{c1}, \eqref{c2}, and \eqref{c3} turns out to be more complex than \eqref{myc2}.
Indeed, as shall be argued shortly, \eqref{myc}, or its equivalent form \eqref{myc2}, can be regarded as the simplest formula attainable for $c_{k,t}^{m,r}$.

Moreover, using \eqref{myc}, one can readily obtain the first few coefficients of highest degree ``by hand''. These are given by
\begin{align*}
c_{k,k+1}^{m,r} & = \frac{m^k}{k+1}, \quad k \geq 0,  \\
c_{k,k}^{m,r} & = m^{k-1} \Big( r - \frac{m}{2} \Big), \quad k \geq 1,  \\
c_{k,k-1}^{m,r} & = \frac{1}{12} k m^{k-2} \big( m^2 -6mr +6r^2 \big), \quad k \geq 2,  \\
c_{k,k-2}^{m,r} & = \frac{1}{12} k(k-1) m^{k-3} r \big( m^2 -3mr +2r^2 \big), \quad k \geq 3, \\
c_{k,k-3}^{m,r} & = \frac{1}{720} k(k-1)(k-2) m^{k-4} \big( m^4 -30m^2 r^2 +60m r^3 -30r^4 \big), \quad k \geq 4,
\end{align*}
etc.

\section{Gould's conjecture}

Setting $x=0$ in \eqref{ber2} yields the following well-known explicit formula for the Bernoulli numbers \cite[Equation 1]{gould}
\begin{equation}\label{gould}
B_{k} = \sum_{j=0}^{k} \sum_{i=0}^{j} \frac{(-1)^{i}}{j+1} \binom{j}{i} i^k, \quad k \geq 0,
\end{equation}
or, equivalently,
\begin{equation}\label{beq}
B_k = \sum_{j=0}^k (-1)^j \frac{j!}{j+1} S_2(k,j),  \quad k \geq 0.
\end{equation}
Regarding the structure of explicit formulas for Bernoulli numbers, it is worth recalling the conjecture made by Gould at the end of his survey paper \cite{gould}:
\begin{itemize}
  \item Gould's conjecture: ``\emph{the writer has seen no formula for $B_k$ which does not require at least two actual summations}.''
\end{itemize}

This statement hints at the recognition that, in fact, there does not exist any elementary formula for $B_k$ involving just one (finite) summation. Hence, in the spirit of Gould's conjecture, we can make the ansatz that the above formula \eqref{gould} for $B_k$ (or its equivalent form \eqref{beq}), as well as the above formula \eqref{ber2} for $B_k(x)$, constitutes, in a suitably defined sense, the simplest explicit formula for the Bernoulli numbers and the Bernoulli polynomials, respectively. Consequently, with this proviso, and noting that
\begin{equation}\label{prov}
c_{k,t}^{m,r} = \frac{m^k}{k+1} \binom{k+1}{t} B_{k+1-t}\left( \frac{r}{m}\right),
\end{equation}
we may conclude that the expression in \eqref{myc}, or its equivalent form \eqref{myc2}, constitutes in turn the simplest formula for $c_{k,t}^{m,r}$ one is ever likely to find.

\section{Concluding remarks}

We end this note with the following three remarks.

\begin{remark}
From \eqref{c1} and \eqref{prov} we can readily obtain an explicit formula for the Bernoulli polynomials involving the Stirling numbers of the first and second kind. Indeed, for $t=1$, \eqref{c1} becomes
\begin{equation*}
c_{k,1}^{m,r} = m^k \sum_{j=1}^{k+1} \sum_{i=j}^{k+1} \frac{j}{i} \left( \frac{r}{m} \right)^{j-1} S_1(i,j) S_2(k,i-1),
\end{equation*}
while, from \eqref{prov}, it follows that
\begin{equation*}
c_{k,1}^{m,r} = m^k B_{k}\left( \frac{r}{m}\right).
\end{equation*}
Therefore, equating the right-hand sides of the last two equations, and replacing the rational $\frac{r}{m}$ by the arbitrary variable $x$, we find that
\begin{equation}\label{beq2}
B_k(x) = \sum_{j=0}^{k} \sum_{i=j}^{k} \frac{j+1}{i+1} S_1(i+1,j+1) S_2(k,i) x^j, \quad k \geq 0,
\end{equation}
which reduces to \eqref{beq} when $x=0$. Let us also observe that, by using the property $B_k(1) = (-1)^k B_k$, from \eqref{beq2} we get the identity
\begin{equation*}
B_k = (-1)^k \sum_{j=0}^{k} \sum_{i=j}^{k} \frac{j+1}{i+1} S_1(i+1,j+1) S_2(k,i), \quad k \geq 0.
\end{equation*}
Moreover, from \eqref{beq2} and $B_k(x) = \sum_{j=0}^{k} \binom{k}{j} B_{k-j} x^j$, it turns out that
\begin{equation*}
\binom{k}{j} B_{k-j} = (j+1) \sum_{i=j}^{k} \frac{1}{i+1} S_1(i+1,j+1) S_2(k,i), \quad 0 \leq j \leq k,
\end{equation*}
which should be compared with the alternative expression \cite[Equation (15.41)]{gould2}
\begin{equation*}
\binom{k}{j} B_{k-j} = k \sum_{i=j}^{k} \frac{1}{i} S_1(i,j) S_2(k-1,i-1), \quad 1 \leq j \leq k.
\end{equation*}
\end{remark}

\begin{remark}
For $t=1$, \eqref{c2}, as well as \eqref{myc2}, becomes
\begin{equation*}
c_{k,1}^{m,r} = \sum_{j=0}^{k} (-1)^j \frac{m^j j!}{j+1} W_{m,r}(k,j),
\end{equation*}
and then
\begin{equation*}
B_{k}\left( \frac{r}{m}\right) = \sum_{j=0}^{k} (-1)^j \frac{j!}{j+1} \frac{W_{m,r}(k,j)}{m^{k-j}},  \quad k \geq 0,
\end{equation*}
which was obtained in \cite{mih} as a consequence of \cite[Theorem 1]{mih}. Moreover, when $m=1$ the last formula reduces to
\begin{equation*}
B_k(r) = \sum_{j=0}^{k} (-1)^j \frac{j!}{j+1} S_{2,r}(k+r,j+r),
\end{equation*}
thus retrieving \cite[Theorem 1]{guo}, where $S_{2,r}(k+r,j+r)$ denotes the $r$-Stirling numbers of the second kind, namely
\begin{equation*}
S_{2,r}(k+r,j+r) = \frac{1}{j!} \sum_{i=0}^{j} (-1)^{j-i} \binom{j}{i} \big( i +r \big)^k .
\end{equation*}
Notice that, by combining the last two equations, we recover \eqref{ber2} with $x =r$.
\end{remark}

\begin{remark}
Clearly, we have that $S_k^{m,r}(1) = r^k$, and then $\sum_{t=1}^{k+1} c_{k,t}^{m,r} = r^k$. Hence, it follows from \eqref{prov} that
\begin{equation*}
\frac{1}{k+1} \sum_{j=0}^{k} \binom{k+1}{j} B_{j}\left( \frac{r}{m}\right) = \left(\frac{r}{m}\right)^k,
\end{equation*}
which generalizes to the well-known relation
\begin{equation*}
\frac{1}{k+1} \sum_{j=0}^{k} \binom{k+1}{j} B_{j}(x) = x^k.
\end{equation*}
\end{remark}


\vspace{.2cm}

\begin{thebibliography}{10}


\bibitem{bazso} A. Bazs\'{o} and I. Mez\H{o}, On the coefficients of power sums of arithmetic progressions, \textit{J. Number Theory}, 153, 117--123, 2015.

\bibitem{cere} J. L. Cereceda, Generalized Stirling numbers and sums of powers of arithmetic progressions, \textit{Internat. J. Math. Ed. Sci. Tech.}, 51(6), 954--966, 2020.

\bibitem{gaut} N. Gauthier, Explicit formula for power sums of an arithmetic sequence, {\it Math. Gaz.}, 91(520), 97--103, 2007.

\bibitem{gould} H. W. Gould, Explicit formulas for Bernoulli numbers, {\it Amer. Math. Monthly}, 79(1), 44--51, 1972.

\bibitem{grif} M. Griffiths, More on sums of powers of an arithmetic progression, {\it Math. Gaz.}, 93(527), 277--279, 2009.

\bibitem{guo} B.-N. Guo, I. Mez\H{o}, and F. Qi, An explicit formula for Bernoulli polynomials in terms of $r$-Stirling numbers of the second kind, {\it Rocky Mountain J. Math.}, 46(6), 1919-1923, 2016.

\bibitem{hirs} M. D. Hirschhorn, Evaluating $\sum_{n=1}^{N} (a +nd)^p$, {\it Math. Gaz.}, 90(517), 114--116, 2006.

\bibitem{mezo} I. Mez\H{o}, A new formula for the Bernoulli polynomials, {\it Results Math.}, 58, 329--335, 2010.

\bibitem{mih} M. Mihoubi and M. Tiachachat, Some applications of the $r$-Whitney numbers, {\it C. R. Acad. Sci. Paris, Ser. I}, 352(12), 965--969, 2014.

\bibitem{gould2} J. Quaintance and H. W. Gould, {\it Combinatorial Identities for Stirling Numbers. The Unpublished Notes of H. W. Gould}. Singapore: World Scientific Publishing, 2016.

\bibitem{ramirez} J. L. Ram\'irez, S. N. Villamar\'in, D. Villamizar, Eulerian numbers associated with arithmetical progressions. {\it The Electronic Journal of Combinatorics}, 25, Article \#P1.48, 2018.

\bibitem{sri} H. M. Srivastava, Some explicit formulas for the Bernoulli and Euler numbers and polynomials, {\it Internat. J. Math. Ed. Sci. Tech.}, 19(1), 79--82, 1988.

\bibitem{todorov} P. G. Todorov, On the theory of the Bernoulli polynomials and numbers, {\it J. Math. Anal. Appl.}, 104(2), 309-350, 1984.

\bibitem{tsao} H. P. Tsao, Explicit polynomial expressions for sums of powers of an arithmetic progression, {\it Math. Gaz.}, 92(523), 87--92, 2008.

\end{thebibliography}
\end{document}